\documentclass{amsart}

\newtheorem{theorem}{Theorem}
\newcommand{\bt}{\begin{theorem}}
\newcommand{\et}{\end{theorem}}

\newtheorem{lemma}{Lemma}
\newcommand{\bl}{\begin{lemma}}
\newcommand{\el}{\end{lemma}}

\newtheorem{corollary}{Corollary}
\newcommand{\bc}{\begin{corollary}}
\newcommand{\ec}{\end{corollary}}

\DeclareMathOperator{\card}{\text{card}}
\DeclareMathOperator{\qand}{\quad\text{and}\quad}
\DeclareMathOperator{\qqand}{\qquad\text{and}\qquad}

\newcommand{\N}{\ensuremath{ \mathbf N }}

\newcommand{\Z}{\ensuremath{\mathbf Z}}

\newcommand{\mba}{\ensuremath{ \mathbf a}}
\newcommand{\mbe}{\ensuremath{ \mathbf e}}
\newcommand{\mbh}{\ensuremath{ \mathbf h}}

\newcommand{\mbA}{\ensuremath{ \mathbf A}}

\newcommand{\beq}{\begin{equation}}
\newcommand{\eeq}{\end{equation}}
\newcommand{\benum}{\begin{enumerate}}
\newcommand{\eenum}{\end{enumerate}}

\usepackage{amsmath,amssymb,amsthm}
\usepackage[all]{xy}
\title{Chromatic sumsets}
\author{Melvyn B. Nathanson}
\address{Lehman College (CUNY), Bronx, NY 10468}
\email{melvyn.nathanson@lehman.cuny.edu}



\date{\today}
\begin{document}
\maketitle

\begin{abstract}
Let $\mbA = (A_1,\ldots, A_q)$  be a $q$-tuple of finite sets of integers. 
Associated to every $q$-tuple of nonnegative integers $\mbh = (h_1,\ldots, h_q)$ 
is the linear form $\mbh\cdot \mbA = h_1 A_1 + \cdots + h_qA_q$.
The set $(\mbh\cdot\mbA)^{(t)}$ consists of all elements of this sumset with at least $t$ representations. 
The structure of the set $(\mbh\cdot\mbA)^{(t)}$ 
is computed for all sufficiently large $h_i$.
\end{abstract}

\section{Coloring the integers}
Let $A$ be a set of integers.  
The $h$-fold sumset $hA$ consists of all integers $n$ that can be represented 
as the sum of $h$ not necessarily distinct elements of $A$.  
The function $r_{A,h}(n)$ counts the number of representations of $n$:  
\[
r_{A,h}(n) = \card\left\{ (a_{j_1},\ldots, a_{j_h} ) \in A^h:  
n = a_{j_1} + \cdots + a_{j_h} \text{ and }  a_{j_1} \leq \cdots \leq  a_{j_h}\right\}. 
\]
Thus, $hA = \{n\in \Z: r_{A,h}(n) \geq 1\}$.  
For every positive integer $t$, let $(hA)^{(t)}$ be the set of integers that have 
at least $t$ such representations:
\[
(hA)^{(t)} = \{n \in \Z: r_{A,h}(n) \geq t\}.  
\]

Let $u$ and $v$ be integers.  Define the \emph{interval of integers} 
\[
[u,v] = \{n \in \Z : u \leq n \leq v\} 
\]
the \emph{dilation}
\[
d\ast A = \{da:a\in A\} 
\]
and the sets
\[
u+A = \{u+a:a\in A\} \qqand v-A = \{v-a:a\in A\}.
\]

Let $A$ be a finite set of integers with 
\[
|A| \geq 2, \quad \min(A) = a_0, \qand d = \gcd\{a-a_0:a\in A\}.
\]
The set 
\[
A_0= \left\{\frac{a-a_0}{d} : a \in A \right\}
\] 
satisfies $\min(A_0) = 0$ and $\gcd(A_0) = 1$.  Moreover,  
\[
(hA)^{(t)} = ha_0 + d\ast (hA_0)^{(t)}.
\]
Thus, the following theorem describes the structure of the sumset $(hA)^{(t)}$ for all finite sets $A$, 
all positive integers $t$, and all sufficiently large $h$.

\bt[Nathanson~\cite{nath1972-7,nath1996bb,nath2021-190}]          \label{chromatic:theorem:MBN}
Let $A$ be a finite set of integers such that    
\[
|A| = k \geq 2, \quad \min(A) = 0, \quad \max(A) = a^*, \qand \gcd(A) = 1.
\]
For every positive integer $t$,  let 
\[
h_t = (k-1)(ta^*-1)a^* + 1.
\]
There are nonnegative integers $c_t$ and $d_t$ 
and finite sets $C_t$ and $D_t$ with 
\[
C_t \subseteq  [0,c_t -2]  \qqand D_t \subseteq [0,d_t -2]
\]
such that 
\[
(hA)^{(t)} = C_t  \cup [c_t,ha^* -d_t] \cup \left(ha^*  - D_t \right) 
\]
 for all $h \geq h_t$.
\et

There is a more subtle additive problem.   
Color the elements of the set $A$ with $q$ colors, 
which we call $\{1,\ldots, q\}$.
Let $A_i$ be the subset of $A$ consisting of all elements $a \in A$ that have color $i$. 
The sets $A_1, \ldots, A_q$ are pairwise disjoint with $A = \bigcup_{i=1}^q A_i$.  
The $q$-tuple $\mbA = (A_1,\ldots, A_q)$ is an ordered partition of $A$.  
Let $\N_0^q$ be the set of  $q$-tuples of nonnegative integers.  
For  
\[
\mbh = (h_1,\ldots, h_q) \in \N_0^q, 
\]  
let
\[
\| \mbh \| = \sum_{i=1}^q h_i = h.
\]

The \emph{chromatic subset}  $\mbh\cdot \mbA$ is the set of all integers in the sumset $hA$ 
that can be represented as the sum of $h$ elements of $A$ 
with exactly $h_i$ elements of color $i$ for all $i \in [1,q]$.  
Thus, $n \in \mbh\cdot \mbA$ if and only if we can write 
\beq                               \label{chromatic:rep-1}
n = \sum_{i=1}^q \sum_{j_i=1}^{h_i} a_{i,j_i}
\eeq
where 
\beq                               \label{chromatic:rep-2}
a_{i,j_i }\in A_i \qquad \text{for all $i \in [1,q]$ and $j_i \in [1,h_i]$.}
\eeq

We refine this problem by allowing elements of $A$ to have more than one color.  
The subset $A_i$ still consists of all elements of $A$ with color $i$, and $A = \bigcup_{i=1}^q A_i$, 
but the sets $A_i$ are not necessarily pairwise disjoint. 
The chromatic sumset $\mbh \cdot \mbA$ consists of all integers $n$ 
that have at least one representation of the form~\eqref{chromatic:rep-1} 
and~\eqref{chromatic:rep-2}. 
Note that 
\[
h_i A_i = \left\{  \sum_{j_i=1}^{h_i} a_{i,j_i}: a_{i,j_i} \in A_i \text{ for all } j_i \in [1,h_i]  \right\} 
\]
and so 
\begin{align*}
\mbh \cdot \mbA & = (h_1,\ldots, h_q)  \cdot (A_1,\ldots, A_q) \\ 
& = h_1A_1 + \cdots + h_qA_q.   
\end{align*}
This is a \emph{homogeneous linear form} in the sets $A_1,\ldots, A_q$.
For every set $B$ of integers, we also have the \emph{inhomogeneous linear form} 
$h_1A_1 + \cdots + h_qA_q + B$.  
Han, Kirfel, and Nathanson~\cite{nath1998-92} 
determined the asymptotic structure of homogeneous and inhomogeneous linear forms 
for all $q$-tuples of finite sets of integers.

The \emph{chromatic representation function} $r_{\mbA,\mbh}(n)$ counts the number 
of colored representations of $n$ of the form~\eqref{chromatic:rep-1} and~\eqref{chromatic:rep-2}, 
that is, the number of $q$-tuples 
\[
\left( \mba_1, \ldots, \mba_q \right) 
\in A_1^{h_1}\times \cdots \times \cdots A_q^{h_q} 
\]
where, for all $i \in [1,q]$, the $h_i$-tuple
 $\mba_i  = (a_{i,j_1}, a_{i,j_2}, \ldots,  a_{i,j_{h_i}} )\in A_i^{h_i}$ satisfies 
\[
a_{i,j_1} \leq a_{i,j_2} \leq \cdots \leq a_{i,j_{h_i}}. 
\]
Let 
\[
(\mbh\cdot \mbA)^{(t)} = \{ n \in \mbh\cdot \mbA : r_{\mbA,\mbh}(n) \geq t\}.
\]
In this paper we determine the structure of the chromatic sumset 
$(\mbh\cdot \mbA)^{(t)}$ for all positive integers $t$ and all sufficiently large vectors $\mbh  \in \N_0^q$. 

Define a partial order on the set $\N_0^q$ as follows.  For vectors 
\[
\mbh_1 = (h_{1,1}, \ldots, h_{q,1}) \in \N_0^q 
\qqand 
\mbh_2 = (h_{1,2}, \ldots, h_{q,2}) \in \N_0^q
\]
let 
\[
\mbh_1 \preceq \mbh_2 \qquad \text{ if $h_{1,i} \leq h_{2,i}$ for all $i \in [1,q]$.}
\]
We also write $\mbh_2 \succeq \mbh_1$ if $\mbh_1 \preceq \mbh_2$.  

Let $\mbh_j = (h_{1,j}, \ldots, h_{q,j}) \in \N_0^q$ for $j \in [ 1, \ell]$.  
Let  
\[
h_i = \sup( h_{i,j}: j \in [ 1, \ell] ) \qquad \text{for $i \in [1,q]$.}
\]
The least upper bound of the vectors $\mbh_1, \ldots, \mbh_{\ell}$  
is the vector  
\[
\sup( \mbh_j : j \in [ 1, \ell] ) = (h_1,\ldots, h_q) \in \N_0^q. 
\]

The $q$-tuple of finite sets of integers $\mbA = (A_1,\ldots, A_q)$ is \emph{normalized} if  
\[
\min(A_i) = 0  \quad\text{for all $i \in [1,q]$}  
\]
and 
\[
\gcd\left( \bigcup_{i=1}^q A_i \right)  = 1.
\]

The case $t=1$ of the following theorem is the result of Han, Kirfel, and Nathanson~\cite{nath1998-92}. 
The case $q=1$ is Theorem~\ref{chromatic:theorem:MBN}.  

\bt                                     \label{chromatic:theorem:structure}
Let $\mbA = (A_1,\ldots, A_q)$ be a normalized $q$-tuple 
of finite sets of integers.  
Let $\max(A_i) = a_i^* \geq 1$ for all $i \in [1,q]$, and
\[
\mba^* = \left( a_1^*,\ldots, a_q^* \right) \in \N_0^q.
\]
For every positive integer $t$, there exist nonnegative integers $c_t$ and $d_t$, 
there exist finite sets of nonnegative integers 
$C_t$ and $D_t$, and there exists a vector $\mbh_t = (h_{t,1},\ldots, h_{t,q} ) \in \N_0^q$ such that, 
if  $\mbh = (h_{1},\ldots, h_q) \in \N_0^q$ and $\mbh \succeq \mbh_t$, then 
\[
\left(\mbh \cdot \mbA\right)^{(t)}  = C_t \cup \left[ c_t, \mbh \cdot \mba^* - d_t  \right] 
\cup \left( \mbh \cdot \mba^* - D_t \right). 
\] 
\et

Let $\mbA' = (A'_1,\ldots, A'_q)$ be a $q$-tuple of nonempty finite sets of integers.  
For all $i \in [1,q]$, let $\min(A'_i) = a'_{i,0}$ and $d = \gcd\left(\bigcup_{i=1}^q (A'_i-a'_{i,0} \right)$.  
Define the sets  
\[
A_i = \frac{1}{d}\ast (A'_i - a'_{i,0}) = \left\{ \frac{a'_i-a'_{i,0}}{d}:a'_i \in A'_i \right\}.  
\]
Equivalently, 
\[
A'_i = d\ast A_i + a'_{i,0}.
\]
The $q$-tuple $\mbA = (A_1,\ldots, A_q)$ is normalized, and, 
for all vectors $\mbh = (h_1,\ldots, h_q) \in \N_0^q$, we have  
\[
\mbh\cdot \mbA' = \sum_{i=1}^q h_iA'_i = d\ast \sum_{i=1}^q h_iA_i + \sum_{i=1}^q h_i  a'_{i,0} 
= d\ast (\mbh\cdot \mbA) + \sum_{i=1}^q h_i  a'_{i,0}.
\]
It follows that Theorem~\ref{chromatic:theorem:structure} describes the asymptotic structure of 
the linear form $\left(\mbh \cdot \mbA'\right)^{(t)}  $ for every $q$-tuple $\mbA'$ of finite sets of integers.

\section{Linear  forms}

\bl                                                        \label{chromatic:lemma:mbh-inclusion}
Let $\mbA = (A_1,\ldots, A_q)$ be a $q$-tuple of sets of integers wih $0 \in A_i$ for all $i \in [1,q]$.
If $\mbh_1, \mbh_2 \in \N_0^q$ and $\mbh_1 \preceq \mbh_2$, 
then $\mbh_1 \cdot \mbA \subseteq \mbh_2 \cdot \mbA$.  
\el

\begin{proof}
Let $\mbh_1 = (h_{1,1}, \ldots, h_{1,q})$ and $\mbh_2 = (h_{2,1}, \ldots, h_{2,q}).$
For all $n = \sum_{i=1}^q n_i \in \mbh_1 \cdot \mbA$ 
with 
\[
n_i = \sum_{j_i=1}^{h_{1,i}}  a_{i,j_i}  \in h_{1,i}A_i
\]
we have  
\[
n_i = \sum_{j_i=1}^{h_{1,i}}  a_{i,j_i}  + (h_{2,i} - h_{1,i})0  \in h_{2,i}A_i
\]
and so $n \in \mbh_2 \cdot \mbA$. 
This completes the proof.  
\end{proof}

\bl                                                                   \label{chromatic:lemma:interval-1}
Let $\mbA = (A_1,\ldots, A_q)$ be a $q$-tuple of sets of integers, and let $A = \bigcup_{i=1}^q A_i$.  
For all vectors $\mbh = (h_1,\ldots, h_q) \in \N_0^q$, 
\benum
\item[(i)]
$
\mbh \cdot \mbA \subseteq \| \mbh \| A.
$
\item[(ii)]
If $A_i \subseteq [u_i ,v_i]$ for all $i \in [1,q]$, then 
\[
\mbh \cdot \mbA \subseteq \left[\sum_{i=1}^q h_i u_i,  \sum_{i=1}^q h_i v_i \right].
\]
\item[(iii)]
If $A_i \subseteq [0,a_i^*]$ for all $i \in [1,q]$ and 
$\mba^* = \left( a_1^*,\ldots, a_q^*\right)$, then 
\[
\mbh \cdot \mbA \subseteq \left[ 0, \mbh \cdot \mba^* \right].
\]
\item[(iv)] 
For every positive integer $t$, 
\[
\bigcup_{\substack{( t_1,\ldots, t_q) \in \N_0^q \\ t_1\cdots t_q \geq t}} 
\left(    \left( h_1 A_1 \right)^{(t_1)} + \cdots +   \left( h_q A_q \right)^{(t_q)}  \right) 
\subseteq \left( \mbh \cdot \mbA \right)^{(t)}.
\]
\eenum
\el

\begin{proof}
(i) For all $i \in [1,q]$, we have $A_i \subseteq A$ and so $h_i A_i \subseteq h_i A$.  Therefore, 
\[
\mbh \cdot \mbA  = \sum_{i=1}^q h_i A_i \subseteq \sum_{i=1}^q (h_i A ) 
= \left(\sum_{i=1}^q h_i  \right) A =  \| \mbh \| A.
\]

(ii) For all $i \in [1,q]$, if $A_i \subseteq [u_i, v_i]$, then $h_iA_i \subseteq [h_i u_i,h_i v_i]$  
and so  
\begin{align*}
\mbh \cdot \mbA 
& = h_1A_1 + \cdots + h_qA_q \\
& \subseteq [h_1u_1, h_1v_1] + \cdots + [h_q u_q, h_qv_q] \\
& = \left[  \sum_{i=1}^q h_i u_i,  \sum_{i=1}^q h_i v_i \right]. 
\end{align*}

(iii)  This follows from (ii) with $u_i = 0$ and $v_i = a_i^*$ for all $i \in [1,q]$,  
and $ \sum_{i=1}^q h_i a_i^* = \mbh \cdot \mba^*$.

(iv) If 
\[
n \in  \left( h_1 A_1 \right)^{(t_1)} + \cdots +   \left( h_q A_q \right)^{(t_q)} 
\]
then for all $i \in [1,q]$  there is an integer $m_i \in \left( h_i A_i \right)^{(t_i)} $  such that 
\[
n = m_1+\cdots + m_q.
\]
The integer $m_i$ has at least $t_i$ representations as a sum of $h_i$ elements of $A_i$, 
and so $n$ has at least $t_1\cdots t_q \geq t$ representations in the sumset 
$h_1A_1 + \cdots + h_q A_q$.  Therefore, $r_{\mbA,\mbh}(n) \geq t$. 
This completes the proof.  
\end{proof}

\bl                                                      \label{chromatic:lemma:intervalSum}
Let $A$ be a finite set of integers with $\min(A)= 0$ and $ \max(A) = a^*$.  
If $c$ and $m$ are integers with $m \geq a^*$, then 
\[
[c,c+ m-1] + A = [c,c+m-1+a^*].
\]
\el

\begin{proof}
Because $\min(A) = 0$ and $\max(A) = a^* $, 
we have 
\[
\min([c,c+ m-1] + A ) = c 
\]
and
\[
 \max([c,c+ m-1] + A ) = c + m-1 + a^*
\]
and so 
\[
[c,c+ m-1] + A \subseteq [c, c + m-1 + a^*].
\]
Because $\{0,a^* \} \subseteq A$ and $c+a^* \leq c+m$, we have 
\begin{align*}
[c,c+ m-1] + A &  \supseteq  [c,c+ m-1] + \{0,a^* \} \\
& = [c,c+ m-1] \cup  [c +a^*,c+m-1+a^*] \\
&  \supseteq  [c,c+ m-1] \cup  [c +m,c+m-1+a^*] \\
& = [c, c+m-1+a^*]. 
\end{align*}
This completes the proof.   
\end{proof}

\bl                                                      \label{chromatic:lemma:translate}
Let $\mbA = (A_1,\ldots, A_q)$ be a $q$-tuple of  sets of integers.  
Let $B$ be a set of integers and let $\mbh$ and $\mbh'$ be vectors in $\N_0^q$. 
If 
\[
S \subseteq (\mbh \cdot \mbA)^{(t)}
\]
 then 
\[
S + B\subseteq (\mbh \cdot \mbA + B)^{(t)}  
\] 
and
\[
S + \mbh' \cdot \mbA\subseteq \left( (\mbh + \mbh') \cdot \mbA \right)^{(t)}.
\] 
\el

\begin{proof}
If $s \in S$ has $t$ distinct representations by the linear form  $\mbh \cdot \mbA$, 
then $s+b$ has $t$ distinct representations by the inhomogeneous linear form $\mbh \cdot \mbA + B$. 

\end{proof}

\section{Reflection and symmetry}
Let $A$ be a finite set of integers with $\min(A) = 0$ and $\max(A) = a^*$.  
Define the \emph{reflection}
\[
\widehat{A} =  \max(A) - A = \{a^* - a: a\in A\}.
\]
The set $A$ is \emph{symmetric} if $A=\widehat{A}$.  

\bl                                  \label{chromatic:lemma:symmetry-1} 
For every finite set $A$ of integers with $0 \in A$ and $\max(A) = a^*$,    
 the reflection $\widehat{A}$ has the following properties:  
\[\min \left( \widehat{A}\right)  = 0, \qquad \max\left( \widehat{A}\right) = a^*,   
\]
\[
 \widehat{\widehat{A}} = A, \qquad \gcd\left( \widehat{A}\right) =  \gcd(A),
 \] 
 and, for all positive integers $h$ and $t$,
\[
 \left(h\widehat{A}\right)^{(t)} = \widehat{ \left(hA\right)^{(t)} }.
\]
\el

\begin{proof}
We have $\min \left( \widehat{A}\right) = a^* - a^* = 0$.   
This implies that 
\[
 \widehat{ \widehat{A} } = \{a^* - \widehat{a}: \widehat{a} \in \widehat{A} \} =  \{a^* - (a^* -a): a\in A\} = A. 
\]
Moreover, $d$ divides $a$ for all $a \in A$ 
if and only if $d$ divides $a^*-a$ for all $a \in A$ 
if and only if $d$ divides $\widehat{a}$ for all $\widehat{a} \in \widehat{A}$, 
and so $ \gcd(A) = \gcd\left( \widehat{A}\right)$.   

If $n \in (h\widehat{A})^{(t)}$, then $r_{\widehat{A},\mbh}(n) \geq t$ and there are $t$ distinct $h$-tuples  
$(\widehat{a}_{i_1,s}, \ldots, \widehat{a}_{i_h,s}) \in \widehat{A}^h$ for $s\in [1,t]$ 
such that 
\[
n = \sum_{j=1}^h \widehat{a}_{i_j,s} = \sum_{j=1}^h \left( a^* - a_{i_j,s}  \right).  
\]
where $a_{i_j,s}  = a^* - \widehat{a}_{i_j,s} \in A$.    Equivalently, 
\[
ha^* - n = \sum_{j=1}^h a_{i_j,s} \in (hA)^{(t)}.  
\]
Therefore, $r_{\mbA,\mbh}(ha^* - n ) \geq t$ and so 
$ha^* - n \in (h A)^{(t)}$ and $n \in \widehat{(hA)^{(t)}}$.
Thus, $ (h\widehat{A})^{(t)} \subseteq \widehat{(hA)^{(t)}}$. 
The proof that $\widehat{(hA)^{(t)}} \subseteq (h\widehat{A})^{(t)}$ is similar.
\end{proof}

\bl                                                      \label{chromatic:lemma:symmetry-2}
Let $\mbA = (A_1,\ldots, A_q)$ be a normalized $q$-tuple 
of finite sets of integers with  
\[
a_i^* = \max A_i  =  \max \widehat{A_i} \qquad \text{for all $i \in [1,q]$.} 
\] 
Let $\widehat{\mbA} = (\widehat{A_1},\ldots, \widehat{A_q})$ and 
\[
\mba^* = \left( a_1^*,\ldots, a_q^* \right) \in \N_0^q.
\]
For all $\mbh \in \N_0^q$ and for all integers $n \in [0,\mbh\cdot \mba^*]$, 
\[
r_{\widehat{\mbA}, \mbh}(n) = r_{\mbA,\mbh} (\mbh\cdot \mba^*-n).
\]
If $c$ and $d$ are nonnegative integers with $c+d \leq \mbh\cdot \mba^*$, then 
\[
[d, \mbh\cdot \mba^* - c] \subseteq \left( \mbh \cdot \widehat{A}\right)^{(t)}
\]
if and only if 
\[
[c, \mbh\cdot \mba^* - d] \subseteq \left( \mbh \cdot A \right)^{(t)}. 
\]
\el

\begin{proof}
Let $\mbh = (h_1,\ldots, h_q) \in \N_0^q$.  
If $r_{\widehat{\mbA}, \mbh}(n) \geq t$, then, for all $i \in [1,q]$, $j_i\in [1,h_i]$, and $s \in [1,t]$, 
there exist integers $\widehat{a}_{i,j_i,s} \in \widehat{A_i}$ such that 
\[
n = \sum_{i=1}^q \sum_{j_i=1}^{h_i} \widehat{a}_{i,j_i,s}
\]
for all $s \in [1,t]$ , and, if $s,s' \in [1,t]$ and $s\neq s'$, 
then there exist ${\ell} \in [1,q]$ and $j_{\ell} \in [1,h_{\ell}] $ 
such that $\widehat{a}_{{\ell},j_{\ell},s} \neq \widehat{a}_{{\ell},j_{\ell},s'}$.  
We have $a_{i,j_i,s}=  a_i^*  - \widehat{a}_{i,j_i,s} \in A_i$, 
and so $a_{{\ell},j_{\ell},s} \neq a_{{\ell},j_{\ell},s'}$.  
From
\[
n = \sum_{i=1}^q \sum_{j_i=1}^{h_i} \left( a_i^* - a_{i,j_i,s} \right) 
= \mbh \cdot \mba^* - \sum_{i=1}^q \sum_{j_i=1}^{h_i} a_{i,j_i,s} 
\]
we obtain 
\[
 \mbh \cdot \mba^* - n = \sum_{i=1}^q \sum_{j_i=1}^{h_i} a_{i,j_i,s} \in (hA)^{(t)}
\]
and so $r_{\mbA,\mbh} (\mbh\cdot \mba^*-n) \geq t$.  
Similarly,  $r_{\mbA,\mbh} (\mbh\cdot \mba^*-n) \geq t$ implies $r_{\widehat{\mbA}, \mbh}(n) \geq t$.      

The observation that $d \leq n \leq \mbh\cdot \mba^* - c$ if and only if 
$c \leq \mbh\cdot \mba^* - n \leq \mbh\cdot \mba^* -d$ completes the proof.  
\end{proof}

\section{Proof of Theorem~\ref{chromatic:theorem:structure}}

\begin{proof}
For all $i \in [1,q]$, let $A_i = \{a_{i,0},a_{i,1},\ldots, a_{i,k_i}\}$, where 
\[
0 = a_{i_0} < a_{i_1} < \cdots < a_{i,k_i} = a_i^*. 
\]
Rearranging the sets $A_i$, we can assume that 
\[
a^* = \max\{a_i^*:i\in [1,q]\} = a_q^* = a_{q,k_q}. 
\] 
Let 
\[
k = \sum_{i=1}^q k_i. 
\]
The divisibility condition $\gcd\left( \bigcup_{i=1}^q A_i \right) = 1$  implies that,  
for every integer $n$ and for all $i \in [1,q]$ and $j_i \in [1,k_i]$, 
there exist  integers $x'_{i,j_i}$ such that 
\[
n = \sum_{i=1}^q \sum_{j_i=1}^{k_i} x'_{i,j_i} a_{i,j_i}. 
\]
For all $s \in [1,t]$, the interval  $[(s-1)a^* , sa^*-1]$ is a complete set of residues modulo $a^*$, 
and so, for all $(i,j_i) \neq (q,k_q)$, there is a unique integer $x_{i,j_i,s}(n)$ such that 
\[
x_{i,j_i,s}(n) \in [(s-1)a^* , sa^*-1]
\]
and 
\[
x_{i,j_i,s}(n) \equiv  x'_{i,j_i} \pmod{a^*}.
\]
It follows that  
\[
n \equiv  \sum_{i=1}^q \sum_{\substack{ j_i =1\\ (i,j_i) \neq (q,k_q)} }^{k_i} x_{i,j_i,s}(n) a_{i,j_i} \pmod{a^*} 
\] 
and so there is an integer $x_{q,k_q,s}(n)$ such that 
\beq                                                               \label{chromatic:n-sum}
n   =  \sum_{i=1}^q \sum_{\substack{ j_i =1\\ (i,j_i) \neq (q,k_q)} }^{k_i} x_{i,j_i,s}(n) a_{i,j_i} + x_{q,k_q,s}(n) a^* 
 =  \sum_{i=1}^q \sum_{ j_i =1 }^{k_i} x_{i,j_i,s}(n) a_{i,j_i}.
\eeq
The inequality  
\begin{align*}
 \sum_{i=1}^q \sum_{\substack{ j_i =1\\ (i,j_i) \neq (q,k_q)} }^{k_i} x_{i,j_i,s}(n) a_{i,j_i} 
& \leq  \sum_{i=1}^q \sum_{\substack{ j_i =1\\ (i,j_i) \neq (q,k_q)} }^{k_i} (sa^*-1) a_i^* \\ 
& \leq  \sum_{i=1}^q k_i (ta^*-1) a^* \\ 
& = k (ta^*-1) a^* 
\end{align*}
implies that if 
\beq                                                       \label{chromatic:big-n}
n \geq k (ta^*-1) a^* 
\eeq 
then 
\[
x_{q,k_q,s}(n) a^*
= n -  \sum_{i=1}^q \sum_{\substack{ j_i =1\\ (i,j_i) \neq (q,k_q)} }^{k_i} x_{i,j_i,s}(n) a_{i,j_i} 
 \geq 0
\]
and so $x_{q,k_q,s}(n)$ is a nonnegative integer.  

Let  $i \in [1,q]$ and $j_i \in [1, k_i]$ with $(i,j_i) \neq (q,k_q)$.  
If $s,s' \in [1,t]$ and $s < s'$, then 
\[
 x_{i,j_i,s}(n) <  sa^* \leq (s'-1)a^* \leq x_{i,j_i}(n,s').  
\]
Thus, if $n$ satisfies inequality~\eqref{chromatic:big-n}, then $n$  
has at least $t$ different representations 
as  a sum of elements of $ \bigcup_{i=1}^q A_i$ with repetitions allowed.  

Let $c_t$ be the smallest integer such that every integer  $n \geq c_t$ 
has at least $t$ different representations 
as  a sum of elements of $ \bigcup_{i=1}^q A_i$ with repetitions allowed.   
Let $C_t$ be the set of all integers $n \in [0,c_t-2]$ such that $n$  has at least $t$ different 
representations as  a sum of elements of $ \bigcup_{i=1}^q A_i$ with repetitions allowed.   
Thus, $C_t \cup [c_t,\infty)$ is the set of all integers that have at least $t$ different 
representations as  a sum of elements of $ \bigcup_{i=1}^q A_i$.  

Let $n \in C_t \cup [c_t,\infty)$, and, for $s \in [1,t]$, let  
\[
n   =  \sum_{i=1}^q \sum_{ j_i =1 }^{k_i} x_{i,j_i,s}(n) a_{i,j_i}
\]  
be $t$ distinct representations of $n$.  
For all $i \in [1,q]$, let  
\[
h_i(n) = \max\left\{ \sum_{ j_i =1 }^{k_i} x_{i,j_i,s}(n) : s \in [1,t] \right\} 
\]
and
\[
\mbh(n) = (h_1(n),\ldots, h_q(n)) \in \N_0^q.
\]
Let
\[
\mbh_0 = \sup\left\{ \mbh(n): n \in C_t \cup [c_t,c_t+ a^* -1]  \right\} 
= \left(h_{0,1}, \ldots, h_{0,q} \right)  \in \N_0^q.
\]
Because $0 \in \bigcap_{i=1}^q A_i$, we have   
\[
r_{\mbA, \mbh_0}(n) \geq r_{\mbA, \mbh(n)}(n) \geq t
\]
 for all $n \in C_t \cup [c_t, c_t + a^* -1]$, and so 
\[
C_t \cup [c_t, c_t+ a^* -1] \subseteq \left( \mbh_0 \cdot \mbA \right)^{(t)} \subseteq \mbh_0\cdot \mbA \subseteq [0, \mbh_0\cdot \mba^* ].
\]
It follows that 
\beq                             \label{chromatic:Ct}
C_t\subseteq \left( \mbh \cdot \mbA \right)^{(t)} 
\eeq
for all vectors $\mbh \in \N_0^q$ with $\mbh \succeq \mbh_0$.

Defining the integer 
\[
d'_t = \mbh_0 \cdot \mba^* - (c_t+ a^* -1) \geq 0   
\] 
gives 
\[
[c_t, \mbh_0 \cdot \mba^*  - d'_t] \subseteq \left( \mbh_0 \cdot \mbA \right)^{(t)}. 
\]
We shall prove that 
\beq                             \label{chromatic:leftGrowth}
[c_t, \mbh \cdot \mba^*  - d'_t] \subseteq \left( \mbh \cdot \mbA \right)^{(t)} 
\eeq
for all vectors $\mbh \in \N_0^q$ with $\mbh \succeq \mbh_0$.  

The proof is by induction on $\ell = \|\mbh - \mbh_0\|$.
If $\ell = 0$, then $\mbh = \mbh_0$ and~\eqref{chromatic:leftGrowth} is true.  
Suppose that~\eqref{chromatic:leftGrowth} is true for all vectors $\mbh \in \N_0^q$ 
with $\mbh \succeq \mbh_0$ and $\| \mbh - \mbh_0\| = \ell$.

Let $\mbe_i = (0,\ldots, 0, 1, 0 , \ldots, 0) \in \N_0^q$ be the vector whose $i$th coordinate is 1 
and whose other coordinates are 0.  
If $\mbh \succeq \mbh_0$ and $\| \mbh - \mbh_0\| = \ell + 1$, then there exists $\ell \in [1,q]$ 
such that $\| \mbh - \mbe_{\ell} -  \mbh_0 \| = \ell$ and $\mbh - \mbe_{\ell} \succeq \mbh_0$. 
The induction hypothesis implies 
\begin{align*}
[c_t,  \mbh \cdot \mba^* -  a_{\ell}^*  - d'_t ]  
=  [c_t,  (\mbh - \mbe_{\ell})\cdot \mba^*  - d'_t] 
 \subseteq  \left( (\mbh - \mbe_{\ell}) \cdot \mbA \right)^{(t)}.
\end{align*}
We have 
\[
\mbh\cdot \mba^* \geq (\mbh_0 +\mbe_{\ell}) \cdot \mba^* =  \mbh_0\cdot \mba^* + a^*_{\ell}   
\]
and so 
\[
\mbh \cdot \mba^* -  a_{\ell}^*  - d'_t  + 1\geq \mbh_0 \cdot \mba^*  - d'_t + 1 
= c_t+ a^* \geq c_t + a^*_{\ell}.   
\]
This inequality  and Lemmas~\ref{chromatic:lemma:intervalSum} and~\ref{chromatic:lemma:translate} imply that 
\begin{align*}
[c_t,  \mbh \cdot \mba^* - d'_t ]  
& = [c_t,  \mbh \cdot \mba^* -  a_{\ell}^*  - d'_t ] \cup 
[\mbh \cdot \mba^* -  a_{\ell}^*  - d'_t +1, \mbh \cdot \mba^*  - d'_t ]    \\ 
& = [c_t,  \mbh \cdot \mba^* -  a_{\ell}^*  - d'_t ] \cup 
[c_t + a_{\ell}^*, \mbh \cdot \mba^*  - d'_t ]    \\ 
& = [c_t,  \mbh \cdot \mba^* -  a_{\ell}^*  - d'_t ]  + \{0,a_{\ell}^* \} \\ 
& = [c_t,  \mbh \cdot \mba^* -  a_{\ell}^*  - d'_t ]  + A_{\ell} \\ 
& \subseteq  \left( (\mbh - \mbe_{\ell}) \cdot \mbA \right)^{(t)} + A_{\ell} \\ 
& \subseteq  \left( (\mbh - \mbe_{\ell}) \cdot \mbA  + A_{\ell} \right)^{(t)} \\ 
& = \left( \mbh  \cdot \mbA \right)^{(t)}.  
\end{align*}
This proves~\eqref{chromatic:leftGrowth} for all $\mbh \succeq \mbh_0$.

Consider the system of reflected sets $\widehat{\mbA} = \left( \widehat{A_1}, \ldots, \widehat{A_q} \right)$.  
Let $d_t$ be the smallest integer such that every integer $n \geq d_t$ has at least $t$ representations 
as the sum of elements of $\bigcup_{i=1}^q \widehat{A_i}$.  
Let $D_t$ be the set of integers $n \in [0,d_t-2]$ that have at least $t$ representations 
as the sum of elements of $\bigcup_{i=1}^q \widehat{A_i}$.  
There exists a vector $\widehat{\mbh}_0 \in \N_0^q$ and a nonnegative integer $c'_t$ such that, 
for all  $\mbh \in \N_0^q$ with $\mbh \succeq \widehat{\mbh}_0$,  
\[
d_t + c'_t \leq \mbh \cdot \mba^* 
\]
and 
\[
D_t \cup [d_t, \mbh \cdot \mba^*  - c'_t] \subseteq \left( \mbh \cdot \widehat{\mbA} \right)^{(t)} .
\]
Applying Lemma~\ref{chromatic:lemma:symmetry-2} gives 
\beq                                                                        \label{chromatic:rightGrowth}
 [c'_t, \mbh \cdot \mba^*  - d_t]  \cup \left( \mbh \cdot \mba^*  - D_t \right) 
 \subseteq \left( \mbh \cdot \mbA \right)^{(t)}. 
\eeq

Let $\mbh_t$ be a vector in $\N_0^q$ such that 
\[
\mbh_t \succeq \sup(\mbh_0, \widehat{\mbh}_0 )
\]
and
\beq                                                                        \label{chromatic:ctdt}
c'_t + d'_t \leq \mbh_t \cdot \mba^*.
\eeq
It follows from~\eqref{chromatic:Ct},~\eqref{chromatic:leftGrowth}, \eqref{chromatic:rightGrowth}, 
and~\eqref{chromatic:ctdt} that 
\[
[c_t, \mbh \cdot \mba^*  - d'_t] \cup  [c'_t, \mbh \cdot \mba^*  - d_t] =  [c_t, \mbh \cdot \mba^*  - d_t]
\]
and, for all $\mbh \succeq \mbh_t$, 
\[
C_t \cup [c_t,  \mbh \cdot \mba^*  - d_t]  \cup \left( \mbh \cdot \mba^*  - D_t \right)  
\subseteq \left( \mbh \cdot \mbA \right)^{(t)}.
\]
The definitions of the sets $C_t$ and $D_t$ imply
\[
C_t \cup [c_t,  \mbh \cdot \mba^*  - d_t]  \cup \left( \mbh \cdot \mba^*  - D_t \right)  
= \left( \mbh \cdot \mbA \right)^{(t)}
\]
for all $\mbh \succeq \mbh_t$.  
This completes the proof, 
\end{proof}

\section{Inhomogeneous linear forms}
Let $\mbA = (A_1,\ldots, A_q)$ be a $q$-tuple of sets of integers, 
and let $B$ be a  set of  integers. 
Associated with each vector $\mbh = (h_1,\ldots, h_q) \in \N_0^q$ 
is the \emph{inhomogeneous linear form of sumsets}\index{linear form}
\[
\mbh \cdot \mbA  + B= h_1A_1 + \cdots + h_qA_q + B.
\]
The \emph{chromatic representation function} $r_{\mbA,\mbh,B}(n)$ counts the number of $(q+1)$-tuples 
\[
\left( \mba_1, \ldots, \mba_i , \ldots, \mba_q , b\right) 
\in A_1^{h_1}\times \cdots \times A_i^{h_i} \times \cdots A_q^{h_q} \times B 
\]
with  $b \in B$ such that, for all $i \in [1,q]$, the $h_i$-tuple
 $\mba_i  = (a_{i,j_1}, a_{i,j_2}, \ldots,  a_{i,j_{h_i}} )\in A_i^{h_i}$ satisfies 
\[
a_{i,j_1} \leq a_{i,j_2} \leq \cdots \leq a_{i,j_{h_i}} 
\]
and 
\[
n = \sum_{i=1}^q \sum_{j_i=1}^{h_i} a_{i,j_i} + b.
\]
We shall determine the structure of the sumset 
\[
(\mbh\cdot \mbA+B)^{(t)} = \{ n \in \mbh\cdot \mbA +B: r_{\mbA,\mbh,B}(n) \geq t\}.
\]

\bt                                   \label{chromatic:theorem:structure-B}
Let $\mbA = (A_1,\ldots, A_q)$ be a normalized $q$-tuple 
of finite sets of integers.  Let  $\max(A_i) = a_i^*$ for all $i \in [1,q]$, and 
\[
\mba^* = \left( a_1^*,\ldots, a_q^* \right) \in \N_0^q.
\]
Let $B$ be a finite set of integers with 
\[
0 = \min(B) < \max(B) = b^*.
\]
For every positive integer $t$, there exist nonnegative integers $c_{t,B}$ and $d_{t,B}$, 
finite sets of nonnegative integers $C_{t,B}$ and $D_{t,B}$, 
and a vector $\mbh_{t,B}  \in \N_0^q$ such that, 
if  $\mbh  \in \N_0^q$ and $\mbh \succeq \mbh_{t,B}$, then 
\[
\left(\mbh \cdot \mbA +B \right)^{(t)}  = C_{t,B} \cup \left[ c_{t,B}, \mbh \cdot \mba^* + b^* - d_{t,B}  \right] 
\cup \left( \mbh \cdot \mba^* + b^* - D_{t,B} \right). 
\] 
\et

\begin{proof}
Let  $\mbh = (h_{1},\ldots, h_q) \in \N_0^q$ and $s \in [1,t]$.  If 
\[
n = \sum_{i=1}^q \sum_{j_i=1}^{h_i} a_{i,j_i,s} 
\]
then 
\[
n+ b = \left( \sum_{i=1}^q \sum_{j_i=1}^{h_i} a_{i,j_i,s} \right) + b
\]
for all $b \in B$.
It follows that if $r_{\mbA,\mbh}(n) \geq t$, then 
$r_{\mbA,\mbh,B}(n+b) \geq t$ for all $b \in B$, and so 
\[
 \left(\mbh \cdot \mbA\right)^{(t)} + B \subseteq  \left(\mbh \cdot \mbA + B\right)^{(t)}.  
\]

By Theorem~\ref{chromatic:theorem:structure}, 
for every positive integer $t$, there exist nonnegative integers $c_t$ and $d_t$  
and a vector $\mbh_t \in \N_0^q$ such that, 
if  $\mbh_t = (h_{1},\ldots, h_q) \in \N_0^q$ and $\mbh \succeq \mbh_t$, then 
\[
 \left[ c_t, \mbh \cdot \mba^* - d_t  \right] \subseteq \left(\mbh \cdot \mbA\right)^{(t)}.
\] 
Choose $\mbh_{t,B}= (h_{B,1},\ldots, h_{B,q})  \succeq \mbh_t$ such that, for all $i \in [1,q]$,  
\[
h_{B,i} \geq c_t \qqand b^* + c_t+ d_t  -1\leq  \mbh_{t,B} \cdot \mba^*.
\]
We have  
\begin{align*}
 \left[ c_t, \mbh_{t,B} \cdot \mba^* + b^*- d_t  \right] 
& =  \left[ c_t, \mbh_{t,B} \cdot \mba^* - d_t  \right] \cup \left[  \mbh_{t,B} \cdot \mba^* - d_t  + 1, \mbh_{t,B} \cdot \mba^* + b^*- d_t  \right] \\
& \subseteq  \ \left[ c_t, \mbh_{t,B} \cdot \mba^* - d_t  \right] \cup \left[ b^* +c_t, \mbh_{t,B} \cdot \mba^* + b^*- d_t  \right] \\
& = \left[ c_t, \mbh_{t,B} \cdot \mba^* - d_t  \right] + \{0,b^*\} \\
&  = \left[ c_t, \mbh_{t,B} \cdot \mba^* - d_t  \right] + B \\
& \subseteq  \left(\mbh_{t,B} \cdot \mbA\right)^{(t)} + B \\
& \subseteq  \left(\mbh_{t,B} \cdot \mbA + B\right)^{(t)}. 
\end{align*}

Let $c_{t,B}$ and $d_{t,B}$ be the smallest nonnegative integers such that 
\[
\left[ c_{t}, \  \mbh_{t,B} \cdot \mba^* + b^*- d_{t}  \right] 
 \subseteq 
  \left[ c_{t,B}, \  \mbh_{t,B} \cdot \mba^* + b^*- d_{t,B}  \right] 
 \subseteq 
 \left(\mbh_{t,B} \cdot \mbA + B\right)^{(t)}. 
\]
Let $C_{t,B}$ be the largest subset of $[0,c_{t,B}-2]$ 
and let $D_{t,B}$ be the largest subset of $[0,d_{t,B}-2]$ 
such that 
\[
\left(\mbh_{t,B} \cdot \mbA + B\right)^{(t)} = 
C_{t,B} \cup   \left[ c_{t,B}, \mbh_{t,B} \cdot \mba^* + b^*- d_{t,B}  \right] 
\cup \left( \mbh_{t,B} \cdot \mba^* + b^* - D_{t,B} \right). 
\]
The proof, by induction on $\|\mbh\|$, that 
\beq
\left(\mbh \cdot \mbA + B\right)^{(t)} = 
C_{t,B} \cup   \left[ c_{t,B}, \mbh \cdot \mba^* + b^*- d_{t,B}  \right] 
\cup \left( \mbh \cdot \mba^* + b^* - D_{t,B} \right) 
\eeq
for all $\mbh \succeq \mbh_{t,B}$, is the same as the proof of Theorem~\ref{chromatic:theorem:structure}. 
\end{proof}

\def\cprime{$'$} \def\cprime{$'$} \def\cprime{$'$}
\providecommand{\bysame}{\leavevmode\hbox to3em{\hrulefill}\thinspace}
\providecommand{\MR}{\relax\ifhmode\unskip\space\fi MR }
\providecommand{\MRhref}[2]{%
  \href{http://www.ams.org/mathscinet-getitem?mr=#1}{#2}
}
\providecommand{\href}[2]{#2}

\end{document}